\documentclass[12pt]{amsart}
\usepackage{a4}
\usepackage{amsthm, amsfonts, amssymb,latexsym,amsmath}
\usepackage{enumerate, color}
\usepackage[english]{babel}%[english]
\usepackage[latin1]{inputenc}

\newtheorem{lemma}{Lemma}

\newtheorem{theorem}{Theorem}

\newcommand{\beq}[1]{\begin{equation}\label{#1}}
\newcommand{\eeq}{\end{equation}}

\title[On sets with few distinct distances]{On sets with few distinct distances}

\author[O. Roche-Newton]{Oliver Roche-Newton}
\address{O. Roche-Newton: 69 Altenburger Stra{\ss}e, Johannes Kepler Universit\"{a}t, Linz, Austria }
\email{o.rochenewton@gmail.com }

\begin{document}

\begin{abstract}
It is widely believed that point sets in the plane which determine few distinct distances must have some special structure. In particular, such sets are believed to be similar to a lattice. This note considers two different ways to quantify this idea. 

Firstly, improving on a result of Hanson \cite{H}, it is proven that if $P= A \times A$ with $A \subset \mathbb R$ and $P$ determines $O(|A|^2)$ distinct distances, then $|A-A|=O\left(|A|^{2-\frac{2}{11}}\right)$. This result gives further evidence that cartesian products which determine few distinct distances have some additive structure.

Secondly, it is shown that if a set $P \subset  \mathbb R^2$ of $N$ points determines $O(N/\sqrt {\log N})$ distinct distances, then there exists a reflection $\mathcal R$ and a set $P' \subset P$ with $|P'| =\Omega ( \log^{3/2} N)$ such that $\mathcal R(P') \subset P$. In other words, sets with few distinct distances have some degree of reflexive symmetry.
\end{abstract}

\maketitle

\section{Introduction}
%\blfootnote{Mathematics Subject Classification 52C10}

Given a set $P$ of $N$ points in $\mathbb R^2$, let $d(P):=\{|p-q|:p,q \in P \}$ be the set of distances determined by $P$.\footnote{Here $|p-q|$ denotes the Euclidean distance between points $p$ and $q$.} A classical and beautiful problem in discrete geometry is the Erd\H{o}s distinct distance conjecture, which that states that $|d(P)| = \Omega (N / \sqrt{\log N})$ for all finite $P \subset \mathbb R^2$. The problem was resolved up to logarithmic factors in a landmark work of Guth and Katz \cite{GK}.

A question that remains wide open concerns the possible structure of point sets which determine few distances. It was suggested by Erd\H{o}s \cite{E} that such a set should have ``lattice structure''. More precisely, he suggested the conjecture that if $|d(P)|$ is minimal then there exists a line which contains $\Omega (N^{1/2})$ points of $P$. The conjecture remains wide open; the current best known result establishes that such a set contains $\Omega (\log N)$ points which are supported on a single line\footnote{See the blog post of Adam Sheffer https://adamsheffer.wordpress.com/2014/10/07/few-distinct-distances-implies-many-points-on-a-line/}.

This paper seeks to find other ways to quantify the qualitative idea that sets with few distinct distances are similar to a lattice. One property of lattice point sets is that they are additively structured. A recent paper of Hanson \cite{H} considered the additive properties of points sets of the form $P=A \times A$ with few distinct distances. The main result in \cite{H} was that
$$|d(A \times A) |=\Omega(|A-A||A|^{1/8}),$$
where
$$A-A:=\{a-b:a,b \in A\}$$
is the \textit{difference set} of $A$. 

In particular, Hanson's result shows that if $|d(A \times A)|=O(|A|^2)$ then $|A-A|=O(|A|^{2-\frac{1}{8}})$. This says that if the set $d(A \times A)$ is small then $A$ must have some degree of additive structure. It appears plausible to conjecture that the exponent $2-\frac{1}{8}$ could be replaced by $1+o(1)$. The first aim of this paper is to improve the result of Hanson and make a step in this direction.

\begin{theorem} \label{thm:main} Let $A \subset \mathbb R$ be a finite set and let $P=A \times A$. Then
$$|d(P)|=\Omega(|A-A|^{11/10}).$$
In particular,
$$|d(P)|=O(|A|^2) \Rightarrow |A-A|=O(|A|^{2-\frac{2}{11}}).$$

\end{theorem}

Although this result is stated in terms of discrete geometry and distance problems, it should perhaps be viewed as a result in additive combinatorics. This is reflected in the notation and tools used for the problem. One can view the Guth-Katz theorem for direct product sets as a sum-product type result. It says that\footnote{See the notation in the next section.} for all finite $A \subset \mathbb R$,
$$|(A-A)^2+(A-A)^2| =\Omega \left( \frac{|A|^2}{\log |A|}\right).$$
Theorem \ref{thm:main} and \cite{H} show that this bound being close to tight implies some additive structure. There is a similarity here with the work of Shkredov \cite{S}, who proved the following inverse sum-product result:
$$|(A-A)(A-A)| =O( |A|^2)  \Rightarrow |A-A| =O \left( |A|^{2-\frac{1}{5}+o(1)} \right).$$

Another feature of lattice point sets is that they are highly symmetric; there exists a reflection which maps $P$ to itself. Therefore, we might expect that any set which determines few distinct distances is in some sense highly symmetric. More precisely, we might expect that such a set $P$ has the property that there is a large subset $P'\subset P$ and some reflection $\mathcal R$ such that $\mathcal R (P') \subset P$. The second aim of this paper is to prove the following result in this direction.

\begin{theorem} \label{thm:main2} Let $P$ be a set of $N$ points in $\mathbb R^2$ such that $|d(P)| \leq N/K$, where $K>1$ is some parameter. Then, there exists a subset $P' \subset P$ with $|P'| =\Omega (K^3)$, and some reflection $\mathcal R$ such that $\mathcal R (P') \subset P$.

In particular, if $|d(P)| = O( N/ \sqrt {\log N} )$, then, there exists a subset $P' \subset P$ with $|P'| =\Omega (\log^{3/2} N)$, and some reflection $\mathcal R$ such that $\mathcal R (P') \subset P$.

\end{theorem}

\section{Notation and Preliminary results}

Throughout this paper, for positive values $X$ and $Y$ the notation $X \gg Y$ is used as a shorthand for $X\geq cY$, for some absolute constant $c>0$.

Similar to the difference set, the \textit{sum set} of $A$ and the \textit{product set} of $A$ are defined respectively as
$$A+A:=\{a+b: a,b \in A\} ,\,\,\,\,\,\,\,\,\, AA=\{ab:a,b \in A\}.$$
The shorthand $2A$ is sometimes used for $A+A$. Similar notation is used for longer combinations of sum and difference set; for example $A+A+A-A-A$ is denoted $3A-2A$. Sets formed by a combination of additive and multiplicative operations on different sets are also considered. For example, if $A,B$ and $C$ are sets of real numbers, then
$$AB+C:=\{ab+c:a\in A, b \in B, c \in C \}.$$
Let $A \subset \mathbb R$ be finite and $\lambda \in \mathbb R$. The set of all dilates of $A$ by $\lambda$ is denoted $\{ \lambda \} A$. That is,
$$\{ \lambda \} A=\{\lambda a : a \in A\}.$$
The curly brackets here are used to distinguish the dilate $\{2\}A$ from the sum set $2A$. Also, $A^2$ denotes the set of all squares of $A$. That is $A^2:=\{a^2: a \in A\}$.

The proof of Theorem \ref{thm:main} is a modification of the argument of Hanson \cite{H}. The key new idea in \cite{H} was the following lemma:

\begin{lemma} \label{Hanson} Let $A\subset \mathbb R$ and let $D=A-A$. Then
$$\{2\}DD \subset 2D^2-2D^2. $$

\end{lemma}

As in \cite{H}, we use the following version of Pl\"{u}nnecke-Ruzsa Theorem (see \cite{P}).

\begin{lemma} \label{Plun} Suppose $A$ is a finite subset of an additive abelian group. Then
$$|mA-nA| \leq \left( \frac{|A+A|}{|A|} \right)^{m+n}|A|.$$

\end{lemma}

We also utilise the following sum-product type result, which follows from the Szemer\'{e}di-Trotter. See Exercise 8.3.3 in Tao-Vu \cite{tv}.

\begin{lemma} \label{AB+C} Let $A,B,C \subset \mathbb R$ be finite sets. Then
$$|AB+C|  \gg (|A||B||C|)^{1/2}.$$
\end{lemma}

To prove Theorem \ref{thm:main2}, we require the following weighted version of the Szemer\'{e}di-Trotter Theorem. The result can be found in the literature, see for example \cite{IKRT}. %but we include the simple proof here for the sake of completeness.

\begin{lemma}\label{wst} Let $P$ be a finite set of points in $\mathbb R^2$ and let $L$ be a set of weighted lines. Each line $l \in L$ is assigned a weight $w(l)$. Let $W_L=\sum_{l \in L} w(l)$ denote the total weight of $L$, and let $w_L=\max_{l\in L} w(l)$ be the maximum weight. Then, the number of weighted incidences $I_w(P,L)$ satisfies
\begin{equation}\label{wind}I_w(P,L):=\sum_{p\in P,\,l\in L: p \in L}w(l)\;\ll \: (w_L)^{1/3} (|P|W_L)^{2/3}  + W_L + w_L|P|.\end{equation}\end{lemma}

%\begin{proof} We dyadically decompose the lines from $L$ according to their weights, and then apply the Szemer\'{e}di-Trotter Theorem for each dyadic class. Recall that the Szemer\'{e}di-Trotter Theorem is the following upper bound for point line incidences:
%$$I(P,L) \ll (|P||L|)^{2/3} +|P| +|L|.$$

%Let $L_j:= \{ l \in L : 2^{j-1} \leq w(l) <2^j \}$. Note that $|L_j| \ll \frac {W_L}{2^j}$. Therefore,

%\begin{align*}
%I_w(P,L)&=\sum_{j \geq 1} \sum_{p \in P, l \in L_j : p \in l} w(l)
%\\& \leq \sum_{j\geq 1} 2^j I(P,L_j)
%\\& \ll \sum_{j \geq 1} 2^j (|P|^{2/3}|L_j|^{2/3} + |P| + |L_j|)
%\\& \ll \left (\sum_{j \geq 1} 2^j |P|^{2/3}|L_j|^{2/3} \right) + |P|w_L+ W_L
%\\& \ll \left (\sum_{j \geq 1} (2^j)^{1/3} |P|^{2/3}W_L^{2/3} \right) + |P|w_L+ W_L
%\\& \ll w_L^{1/3} |P|^{2/3}W_L^{2/3} + |P|w_L+ W_L.
%\end{align*} 

%\end{proof}

\section{Proof of Theorem \ref{thm:main}}

Let $D=A-A$ and note that $d(A \times A)=D^2+D^2$. By Lemma \ref{Hanson} and Lemma \ref{Plun}
$$|\{2\}DD+D^2|\leq |3D^2-2D^2| \leq \left( \frac{|D^2+D^2|}{|D^2|} \right )^5 |D^2|.$$
By Lemma \ref{AB+C}, 
$$|\{2\}DD+D^2| \gg |D|^{3/2}.$$
Combining these two estimates, we have
$$|d(A \times A)|=|D^2+D^2| \gg |D|^{11/10}$$
as required.

\section{Proof of Theorem \ref{thm:main2}}

The proof makes use of some observations from a recent paper of Lund, Sheffer and de Zeeuw \cite{LSdZ}, which considered structural properties of point sets which determine few distinct distances via studying the perpendicular bisectors determined by $P$.

We will double count the set of (ordered) isosceles triangles determined by $P$. That is, the proof proceeds by comparing an upper and lower bound for the quantity
$$T:=\{(p,q,s) \in P \times P \times P : |p-s|=|q-s|, p \neq q \}.$$
For the upper bound we use Lemma \ref{wst}. For two distinct points $p,q \in \mathbb R^2$, let $B(p,q)$ denote their perpendicular bisector. Let $L$ be the multiset of perpendicular bisectors determined by $P$. For $l \in L$, its weight is the number of (ordered) pairs of points from $P$ that determine $l$ as a bisector; that is,
$$w(l):=\{(p,q) \in P \times P : B(p,q)=l \}.$$
Note that $W_L = |P|^2-|P| < |P|^2$. Note also that $T=I_w(P,L)$. Indeed $(p,q,s) \in T$ if and only if $s \in B(p,q)$. Therefore, it follows from Lemma \ref{wst} that 
\begin{equation}
|T| \ll w^{1/3}_LN^2 +N^2 \ll w^{1/3}_LN^2.
\label{Tbound}
\end{equation}

On the other hand, if we denote by $C(s,r)$ the circle with radius\footnote{The possibility that $r=0$ is not excluded.} $r$ and centre $s$, then
\begin{align*}
|T|&= \sum_{s\in P} \sum _{r \in d(P)} 2{|C(s,r) \cap P| \choose 2} 
\\& \gg  \sum_{s \in P} \sum_{r \in d(P)} |C(s,r) \cap P|^2-\sum_{s \in P, r \in d(P) : |C(s,r) \cap P| \leq 1} 1 
\\& \gg \sum_{s \in P} \sum_{r \in d(P)} |C(s,r) \cap P|^2-\frac{N^2}{K} .
\end{align*}
Combining this information with \eqref{Tbound}, we deduce that
\begin{equation}
\sum_{s \in P} \sum_{r \in d(P)} |C(s,r) \cap P|^2 = \ll w^{1/3}_LN^2+\frac{N^2}{K}  \ll w^{1/3}_LN^2.
\label{Tbound2}
\end{equation}
Note also that
$$\sum_{s \in P} \sum_{r \in d(P)} |C(s,r) \cap P| = N^2.$$
Therefore, by the Cauchy-Schwarz inequality and \eqref{Tbound2}, we have
\begin{align*}N^4 &= \left(\sum_{s \in P} \sum_{r \in d(P)} |C(s,r) \cap P|\right)^2 
\\&\leq N|d(P)|\sum_{s \in P} \sum_{r \in d(P)} |C(s,r) \cap P|^2
\\& \ll \frac{N^4w_L^{1/3}}{K}.
\end{align*}
This tells us that $w_L =\Omega(K^3)$.  This completes the proof, since there is some perpendicular bisector $l$ such that $l=B(p_i,q_i)$ for $i=1,\dots,k$, $k=\Omega(K^3)$ and with the $p_i$ all distinct. Therefore we can take $P'=\{p_i:1 \leq i \leq k\}$ and $\mathcal R$ to be reflection in the line $l$, and observe that
$$\mathcal R(P')=\{q_i:1 \leq i \leq k\} \subset P.$$

\section{Acknowledgements} The author was supported by the Austrian Science Fund FWF Project F5511-N26,
which is part of the Special Research Program "Quasi-Monte Carlo Methods: Theory and Applications". With thanks to Brandon Hanson and Adam Sheffer for helpful discussions. Much of this work originated from discussions at the IPAM reunion conference for the program ``Algebraic Techniques for Combinatorial and Computational Geometry''. IPAM is funded by the NSF.

%\subsection*{Acknowledgements}J.~Cilleruelo  was supported by grants MTM 2011-22851 of MICINN and ICMAT Severo
%Ochoa project SEV-2011-0087. Oliver Roche-Newton was supported by EPSRC Doctoral Prize Scheme (Grant Ref:  EP/K503125/1) and by the Austrian Science Fund (FWF): Project F5511-N26, which is part of the Special Research Program ``Quasi-Monte Carlo Methods: Theory and Applications.


\begin{thebibliography}{99}

\bibitem{E} P. Erd\H{o}s, {\it On some metric and combinatorial geometric problems}.  Discrete Math. {\bf 60} (1986), 147-152.

\bibitem{GK} L. Guth,  N. H. Katz, {\it On the Erd\H os distinct distance problem in the plane}.   Ann. of Math. (2) {\bf 181} (2015), no. 1, 155--190.

\bibitem{H} B. Hanson, `The Additive Structure of Cartesian Products Spanning Few Distinct Distances',  Preprint arXiv:1607.03442.

\bibitem{IKRT} A. Iosevich, S. Konyagin, M. Rudnev, M, V. Ten, `Combinatorial complexity of convex sequences', Discrete Comput. Geom. {\bf 35} (2006), no. 1, 143--158.

\bibitem{LSdZ} B. Lund, A. Sheffer and F. de Zeeuw, `Bisector energy and few distinct distances',  Preprint arXiv:1411.6868.

\bibitem{P} G. Petridis, `New proofs of Pl\"{u}nnecke-type estimates for product sets in groups', Combinatorica {\bf 32} (2012), no. 6, 721-733.

\bibitem{S} I. Shkredov, `Difference sets are not multiplicatively closed',  Preprint arXiv:1602.02360.

\bibitem{tv} T. Tao, V. Vu. 'Additive combinatorics' \textit{Cambridge University Press} (2006).








\end{thebibliography}
\end{document}